\documentclass[Journal,9pt]{IEEEtran}

\IEEEoverridecommandlockouts                              
\usepackage{graphicx}
\usepackage{multicol}
\usepackage{graphicx}
\usepackage{color}
\usepackage[ansinew]{inputenc}         
    \usepackage{balance}                

\usepackage[all]{xy}
\usepackage{cite}
\usepackage{multicol}
\usepackage{subfigure}
\usepackage{amsthm}
\usepackage{mathrsfs}
\usepackage{amsmath,amssymb,epsf,epsfig,times}
\usepackage{amssymb}
\usepackage{graphicx}
\usepackage{latexsym,amsmath,amsfonts}
\usepackage{mathrsfs}
\usepackage{amscd}
\usepackage{color}

\def\QED{\mbox{\rule[0pt]{1.5ex}{1.5ex}}}
\def\endproof{\hspace*{\fill}~\QED\par\endtrivlist\unskip}

\newtheorem{theorem}{Theorem}

\newtheorem{remark}{Remark}
\newtheorem{lemma}{Lemma}
\newtheorem{assumption}{Assumption}

\title{\LARGE \bf Distributed velocity-constrained consensus of discrete-time multi-agent
systems with nonconvex constraints, switching topologies, and delays\thanks{This work was supported by the National Science Foundation under Grant CMMI-1537729, the National Natural Science Foundation of China (61203080,61573082,61528301,61120106010), and the State Key Laboratory of Intelligent Control and Decision of Complex Systems of Beijing Institute of Technology.}}

\author{Peng Lin$^{a}$, Wei Ren$^{b}$, Huijun Gao$^c$
\thanks{$^{a}$Peng Lin is with the School of Electronic and Information Engineering, Xi'an Jiaotong University, Xi'an, China.
        {\tt\small lin$\_$peng0103@sohu.com}}%
\thanks{$^{b}$Wei Ren is with the Department of Electrical
Engineering, University of California, Riverside, USA.
        {\tt\small ren@ee.ucr.edu}}%
        \thanks{$^{c}$Huijun Gao is with the Research Institute of Intelligent Control and Systems, Harbin Institute of Technology.
        {\tt\small huijungao@hit.edu.cn}}
}

\begin{document}

\renewcommand{\baselinestretch}{0.95}

\maketitle \thispagestyle{empty} \pagestyle{empty}

\begin{abstract}
In this paper, a distributed velocity-constrained consensus problem is studied
for discrete-time multi-agent systems, where each agent's velocity is constrained to lie in a nonconvex set. A distributed constrained control algorithm
is proposed to enable all agents to converge to a common point using only local information. {The gains of the algorithm for all agents need not to be the same or predesigned and
can be adjusted by each agent itself based on its own and neighbors' information.}
  It is shown that the algorithm is robust to arbitrarily bounded communication delays and arbitrarily switching communication
graphs provided that
 the union of the graphs has
directed spanning trees among each certain time interval. The analysis approach is based on
multiple novel model transformations, proper control parameter selections, boundedness analysis of state-dependent stochastic matrices{\footnote{Throughout this paper, when referring to a stochastic matrix, it means a row stochastic matrix.}},
exploitation of the convexity of stochastic matrices, and
the joint connectivity of the communication graphs. 
 Numerical examples are included to illustrate
the theoretical results.

\vspace{0.1cm}

\end{abstract}

\section{Introduction}
In recent years, consensus problems in multi-agent systems have
received a great deal of attention because of its important
applications including formation control, attitude alignment of
clusters of satellites, and flocking \cite{saber, xiao, lin74, hong, xie, ren1, yangmeng, Lixiangwei, Mengzhaolin, CaoRen, angelia, LinRen12-a, srivast, NedicOzdaglar10, plinwrenysong}.
Most of the
existing results concentrate on the ideal case where the state or input of each
agent has no constraints. In some practical situations, the state or input of
each agent is usually constrained to lie in a certain set, e.g., the saturation and dead zone
of the velocity of physical vehicles.

Research on consensus problems with
state or input constraints can be found in \cite{ren1,yangmeng,Lixiangwei,Mengzhaolin, CaoRen,angelia,LinRen12-a,srivast}. For example, article
\cite{ren1} introduced hyperbolic tangent functions to a consensus algorithm for continuous-time double-integrator multi-agent systems with a fixed undirected
topology where the maximum amplitude of the control input of each agent is upper bounded. Also, from the view point of saturation control, articles
 \cite{yangmeng,Lixiangwei,Mengzhaolin,CaoRen} studied constrained control problems by a Lyapunov approach and showed that consensus can be achieved asymptotically or in finite time.
However, in \cite{ren1,yangmeng,Lixiangwei,Mengzhaolin,CaoRen}, it is assumed that each agent has continuous-time dynamics, the input constraint set of each agent is a hypercube and the communication graph
is undirected. From the view point of projection control,
article \cite{angelia}
proposed a projection algorithm for discrete-time multi-agent systems with switching topologies, where each agent is assumed to remain in a convex set.
 Founded on \cite{angelia}, article \cite{LinRen12-a} took the communication delays into account and showed that the projection algorithm is robust to
 arbitrarily bounded communication delays,
 while article \cite{srivast} studied the projection algorithm in a random environment
 and introduced a step size sequence for the consensus stability of the systems.
{However, in \cite{angelia,LinRen12-a,srivast}, it is assumed that the states of the agents are constrained in certain convex sets. When
more general constraint sets are taken into account, the results in \cite{angelia,LinRen12-a,srivast} cannot be directly applied
due to the {loss of the convexity} of the constraint sets.}
%
%


In this paper, our objective is to solve the velocity-constrained consensus
problem for discrete-time multi-agent systems with switching
topologies and nonuniform communication delays. In contrast to \cite{ren1,yangmeng,Lixiangwei,Mengzhaolin, CaoRen}, where the constraint set of each agent is a hypercube,
here each agent' velocity is constrained to lie in a nonconvex set.
The communication graph considered is directed coupled with arbitrarily bounded communication delays and
can be arbitrarily switching as long as the union of the graphs has
directed spanning trees among each certain time interval. To solve the velocity-constrained consensus
problem in this setting, 
 a distributed control algorithm
is proposed by applying a constrained control scheme using only local information. {The gains of the algorithm for all agents need not to be the same or predesigned and
can be adjusted by each agent itself based on its own and neighbors' information,} which distinguishes it from the existing works on double-integrator consensus \cite{xie,hong,lin74}, where the feedback gains
are uniform for all agents.
Owing to the coexistence of the
coupling of the position and velocity states and a velocity delay during the updating process of the position states, the nonlinearity caused
by the nonconvex constraints  would further lead to a stronger nonlinearity on the position states. Both nonlinearities are greatly different from those in  \cite{angelia,LinRen12-a,srivast} and the approaches there cannot be directly applied.
Our analysis approach is to introduce multiple novel model transformations and select proper control parameters to transform the original system into an equivalent system whose system matrix is a state-dependent stochastic matrix. The state-dependent stochastic matrix has two features: one is that the nonzero entries are from an infinite set and the nonzero entries might not be uniformly lower bounded by a positive constant, and the other is that
the stochastic matrix has zero diagonal entries. The coexistence of these two factors poses significant challenges. Most of the existing results
on delay-related consensus require the number of possible nonzero entries to be finite (e.g., \cite{xiao,lin74,LinRen12-a}) and existing approaches
based on the results in \cite{wolf} require the stochastic matrices to have positive
diagonal entries and their nonzero entries to be uniformly lower bounded by a positive constant. Though the results of \cite{NedicOzdaglar10}
allow for an infinite number of edge weights and zero diagonal entries, the union of the communication graphs among each certain time interval is assumed to be strongly connected and each nonzero entry of the stochastic matrices is assumed to be uniformly lower bounded by a positive constant. As a result, the existing results cannot be directly applied to deal with the problem studied in this paper. To study the consensus stability of the equivalent system, we construct an auxiliary matrix each entry of which is no larger than that of the transition matrix of the equivalent system. By analyzing the graph connectivity, we show that the auxiliary matrix  and hence the transition matrix of the equivalent system  have at least one column with all positive entries over a certain time interval. Then, we use the convexity of a stochastic matrix to study the convergence of the transition matrix and show that all its rows tend to the same exponentially as time evolves.


\section{Notations and Preliminaries}
In this section, we introduce some notations and preliminary results
on graph theory and nonnegative matrices (referring to \cite{s10} and \cite{s11}).

{\bf Notations.} $\mathbb{R}^m$  denotes the set of all $m$ dimensional
real column
vectors; $I_m$ denotes the $m$ dimensional unit matrix; $\mathbb{Z}$ denotes the set of all integers; 
$\otimes$ denotes the Kronecker product; $x^T$ denotes the
transpose of a vector $x$; $\mathrm{diag}\{A_1,\cdots,A_q\}$ is a block diagonal matrix with its diagonal blocks equal to the matrices $A_i$, $i=1,\cdots,q$;
{{$\mathrm{\overline{diag}}\{A\}$ denotes a diagonal matrix whose diagonal entries are equal to those of $A$ correspondingly}};
$\inf_{x\in X}x$ denotes the infimum of $x$ in the set $X$; $\prod_{i=s}^kA_i=A_k\cdots A_s$ denotes the product of the matrices $A_k,\cdots,A_s$; $\textbf{1}$
represents
a column vector of all ones with a compatible dimension;
 $\|x\|$ denotes
the standard Euclidean norm of a vector $x$; 0 denotes a zero vector or zero
matrix with an appropriate dimension; $x_i$
denotes the $i$th entry of a vector $x$; and $\lfloor A\rfloor_{ij}$
denotes the $ij$th entry of a matrix $A$.

Let $\mathcal{G}(\mathcal{V},\mathcal{E})$ be a directed
graph of order $n$, where $\mathcal{V}=\{1,\cdots,n\}$ is the set of
nodes, and $\mathcal{E}\subseteq\mathcal{V}\times \mathcal{V}$ is the
set of ordered edges. An edge of $\mathcal{G}$, denoted by $(j,i)$, denotes that agent $i$ can obtain information from agent $j$ but not necessarily vice versa.
Then
the set of neighbors of node $i$ is denoted by $\mathcal{N}_i=\{j\in
\mathcal{V}:(j,i)\in \mathcal{E}\}$. {{The edge weight of each edge $(j,i)$ is defined such that  $a_{ij}>0$ if $(j,i) \in \mathcal{E}$ and $a_{ij}=0$ otherwise.}} The Laplacian of the directed graph $\mathcal{G}$, denoted by $L$,
is defined as $\lfloor {L}\rfloor_{ii} =\sum_{j=1}^na_{ij}$ and $\lfloor {L}\rfloor_{ij}=-a_{ij}$ for all $i\neq j$. The union of a collection of graphs is a
graph whose node and edge sets are the unions of the node and edge sets of the graphs in
the collection.
  A directed path is a sequence of ordered edges of the
form $({i_1},{i_2}),({i_2},{i_3}),\cdots ,$ where ${i_j}\in
\mathcal{V}$ in a directed graph. A directed graph is
strongly connected if there is a directed path from every node to
every other node.  A directed graph has a directed spanning tree, if there exists at
least one node that has directed paths to all other nodes. The node  that has directed paths to all other nodes
is called the root of the directed spanning tree.

Given $C=[c_{ij}]\in \mathbb{R}^{n\times r}$,
$C$ is nonnegative ($C\geq0$) if all its elements $c_{ij}$ are
nonnegative, and  $C$ is positive ($C>0$) if all its
elements $c_{ij}$ are positive. Furthermore, $C\geq D$ if $C-D\geq0$,
and $C>D$ if $C-D>0$. If a nonnegative matrix $C\in
\mathbb{R}^{n\times n}$ satisfies $C\textbf{1}=\textbf{1}$, then it
is {stochastic}.
\section{Model and Problem Statement}
 Consider a multi-agent system consisting
of $n$ agents with discrete-time dynamics. Each agent is regarded as
a node in a swiching directed graph $\mathcal{G}(kT)$, where $k$ is the discrete time index and $T$ is the
sampling period. The Laplacian of the directed graph $\mathcal{G}(kT)$ is denoted by ${L}(kT)$. Each agent updates its current state based
upon the information received from its neighbors, denoted by
$\mathcal{N}_{i}(kT)$. In \cite{angelia,LinRen12-a,srivast}, agents with single-integrator discrete-time dynamics were studied for constrained consensus. In reality, agents usually have double-integrator dynamics and  corresponding constrained consensus has been studied in \cite{ren1,yangmeng,Lixiangwei,Mengzhaolin, CaoRen}, but the constraint sets are limited to hyperplanes. To this end, we study the case of nonconvex constraint sets for the agents with double-integrator discrete-time dynamics which have the following form:
\begin{eqnarray}\label{eq11}\begin{array}{lll}
{x}_i((k+1)T)&=&x_i(kT)+v_i(kT)T\\
{v}_i((k+1)T)&=&u_i(kT)
\end{array}\end{eqnarray}
where $x_i(kT)\in \mathbb{R}^r$ and $v_i(kT)\in \mathbb{R}^r$ for some positive integer $r$ are
the position and velocity states of agent $i$ and $u_i(kT)\in \mathbb{R}^r$ is
the control input. To simplify the notations, we replace all
``$(kT)$" by ``$(k)$". {It is assumed that the initial conditions of $x_i(k)$ and
$v_i(k)$ for all $k\leq 0$ and all $i$ satisfy the dynamics of (\ref{eq11}), and the velocity state of each agent
$v_i(k)$ is constrained to lie in a nonempty constraint set
$V_i\subseteq \mathbb{R}^r$ known only to agent $i$.}

\begin{figure}
 \centering
\includegraphics[width=1.7in]{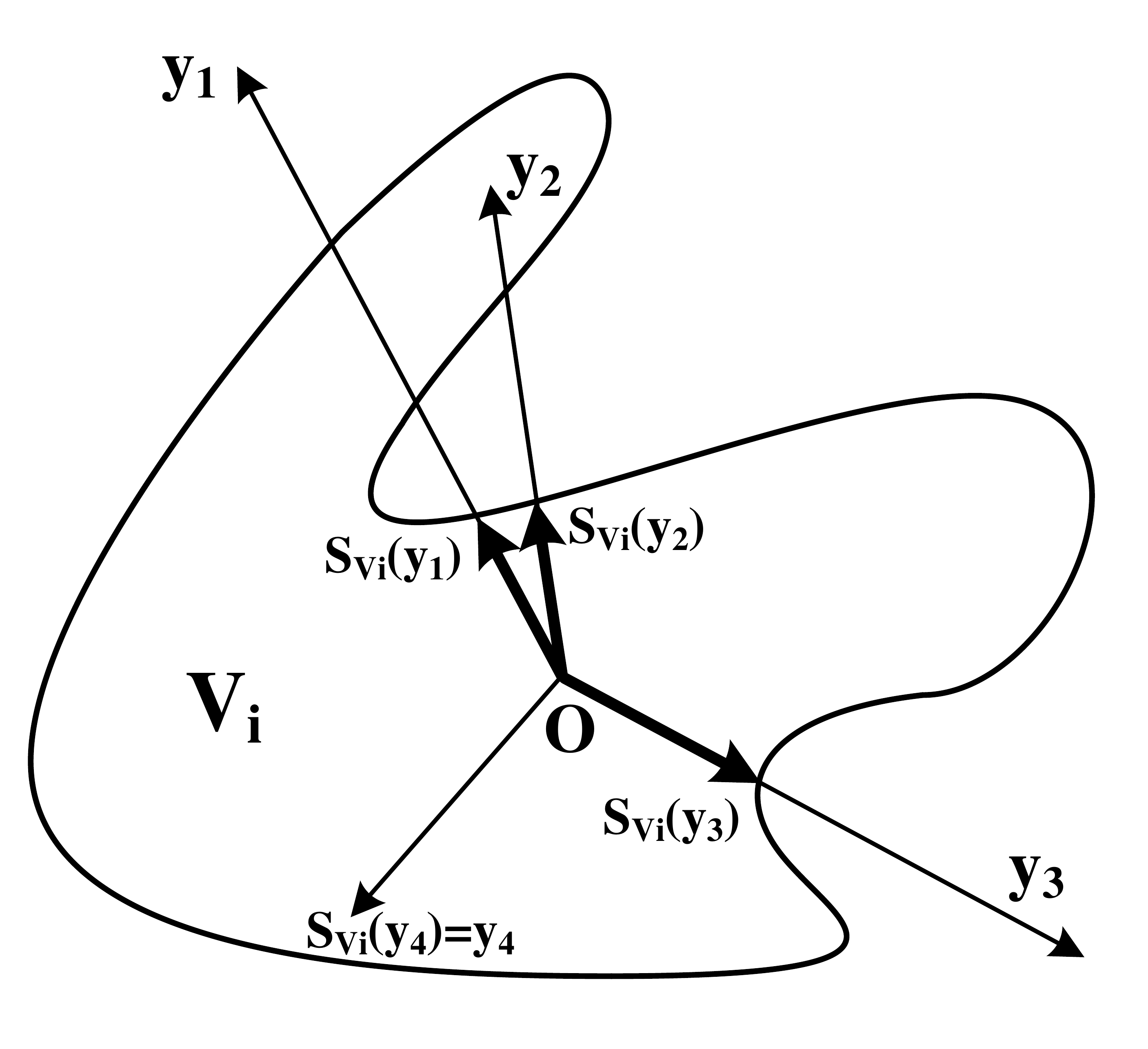} \\
\vspace{-0.4cm}
\caption{Four examples of the constraint operator.}
\label{fig:a1}
\end{figure}

 {Due to the different constraints of each agent's driving forces in different directions, the velocities of the agents, e.g., quadrotors, might not lie in convex sets. Hence we make
the following assumption for $V_i$:}
 \begin{assumption}\label{ass2}{\rm Let $V_i\subseteq \mathbb{R}^r, i=1,\cdots,n$, be nonempty bounded closed
 sets such that $0\in V_i$, $\max_{x\in V_i}\|\mathrm{S}_{V_i}(x)\|=\bar{\rho}_i>0$ and $\inf_{x\notin V_i}\|\mathrm{S}_{V_i}(x)\|=\underline{\rho}_i>0$ for all $i$, where $\bar{\rho}_i$ and
 $\underline{\rho}_i$ are two positive
constants, and $\mathrm{S}_{V_i}(\cdot)$ is a constraint operator such that $\mathrm{S}_{V_i}(0)=0$ and
$\mathrm{S}_{V_i}(x)=\frac{\displaystyle x}{\displaystyle\|x\|}
\max_{{0\leq\beta\leq\|x\|}}{\Big\{}\beta{\big|}\frac{\displaystyle \alpha \beta x}{\displaystyle\|x\|}\in V_i, \forall 0\leq \alpha\leq 1{\Big\}}$ when $x\neq 0$.
}\end{assumption}
{The operator $\mathrm{S}_{V_i}(x)$ means to find the vector with the largest magnitude such that $\mathrm{S}_{V_i}(x)$ has the same direction as $x$, $\|\mathrm{S}_{V_i}(x)\|\leq\|x\|$ and $\alpha\mathrm{S}_{V_i}(x)\in V_i$ for all $0\leq \alpha\leq 1$. (See Fig. \ref{fig:a1} for illustrations.) It should be noted here that we do not impose any {{convexity assumption}} on each $V_i$.
The maximum $\max_{x\in V_i}\|\mathrm{S}_{V_i}(x)\|=\bar{\rho}_i>0$ means that the distance from any point in $V_i$ to the origin is upper bounded. That is, the velocities of all agents cannot be arbitrarily large. The infimum
 $\inf_{x\notin V_i}\|\mathrm{S}_{V_i}(x)\|=\underline{\rho}_i>0$ means that the distance from any point outside $V_i$ to the origin is lower bounded by a positive constant. That is, each agent can move in any direction. Future work could be directed to the more general case where each agent might not be able to move toward certain directions.}

 {{Under the constraint that $v_i(k)\in V_i$ for all $i,k$, our objective is to design an algorithm for all agents to cooperatively reach a consensus on their position states at some vector, denoted by $\bar{x}\in \mathbb{R}^r$, as $k\rightarrow+\infty$, i.e., $\lim_{k\rightarrow+\infty}x_i(k)=\bar{x}$ for all $i$. From the dynamics of (\ref{eq11}), the limit $\lim_{k\rightarrow+\infty}x_i(k)=\bar{x}$ for all $i$ means that $\lim_{k\rightarrow+\infty}v_i(k)=0$ for all $i$. That is, the velocity of each agent would actually converge to zero as $k\rightarrow+\infty$.}}

\section{Main Results}
In this section, we study the velocity-constrained consensus problem for
discrete-time multi-agent systems with switching topologies and communication delays. {Motivated by the algorithms in \cite{lin74,hong,xie} dealing with the case without constraints, we
propose the control algorithm:}
\begin{eqnarray}\label{eq526}\begin{array}{lll}u_i(k)=\mathrm{S}_{V_i}[v_i(k)-p_i(k)v_i(k)T+\pi_i(k)]\end{array}\end{eqnarray} for all $k\geq 0$,
where $v_i(0)=S_{V_i}(v_i(0))$, $p_i(k)>0$ is the feedback damping gain of agent $i$, $\pi_i(k)=\sum\limits_{j\in \mathcal{N}_i(k)}a_{ij}(k)(x_j(k-\tau_{ij}(k))-x_i(k))T$, $0\leq\tau_{i{j}}(k)\in \mathbb{Z}$ is the communication delay from
agent $j$ to agent $i$,
and $a_{ij}(k)$ denotes the edge weight of the edge $(j,i)$
($a_{ij}(k)>0$ for all $j\in \mathcal{N}_i(k)$). It is assumed that all $\tau_{i{j}}(k)$ are upper bounded, i.e., $\tau_{i{j}}(k)\leq M$ for  some constant $M>0$. {When there are no constraints, the algorithm (\ref{eq526}) would have the form of the algorithms introduced in \cite{lin74,hong,xie}. If one agent
receives multiple pieces of the state information from another agent at time $k$,
the latest piece would be used and all others  dropped.
 Here, it is assumed that $a_{ij}(k)\geq \mu_c$ for some
positive number $\mu_c$ when $a_{ij}(k)>0$. The constraint operator is used to ensure the velocity of each agent to be lying in its corresponding constraint set,
and the algorithm parameters of all agents, $p_i(k)$, need not remain the same and
will be shown to be able to be adjusted by each agent itself based on the scaling factors of the constraint operator and
 the parameters of the previous time instant.}  

 \begin{remark}{\rm In \cite{angelia,LinRen12-a,srivast}, the projection operator is used to guarantee all agents with single-integrator discrete-time dynamics remain in their constraint sets. Different from \cite{angelia,LinRen12-a,srivast}, the system (\ref{eq11}) takes the double-integrator form and there is a velocity delay inherent in the dynamics during the updating process of the position states. If the projection operator were used in (\ref{eq526}), due to the coexistence of the nonconvexity of the constraint sets, the
coupling of the position and velocity states, and the velocity delay, the nonlinearity caused by the projection operator would be hard to be measured or estimated and thus the system might become too complicated to analyze. Hence we do not adopt the projection operator in (\ref{eq526}).}\end{remark}

Define \begin{eqnarray}\label{ei10}e_i(k)=\frac{\|\mathrm{S}_{V_i}[v_i(k)-p_i(k)v_i(k)T+\pi_i(k)]\|}{\|v_i(k)-p_i(k)v_i(k)T+\pi_i(k)\|}\end{eqnarray} for all $k\geq0$. In particular, when $v_i(k)-p_i(k)v_i(k)T+\pi_i(k)=0$, we define $e_i(k)=1$. Clearly, $0<e_i(k)\leq 1$. Let $b_i(k)=\frac{1-e_i(k)(1-p_i(k)T)}{T}$. We make the following assumption.

 \begin{assumption}\label{ass3}{\rm Suppose that $\frac{1}{T}>p_i(k+1)\geq b_i(k)>0$ for all $k\geq0$ and all $i$, and there exist a constant $d_{i}>0$ such that $p^2_i(k)> 4 d_{i}\geq 4 \lfloor {L}({k})\rfloor_{ii}$ for all $i$ and all $k\geq0$.}\end{assumption}

To illustrate, we show how to select $p_i(k)$ in a distributed manner to guarantee Assumption \ref{ass3} in three steps:

1. Select $p_i(0)$ such that $0<p_i(0)T<1$. Then calculate $b_i(0)$ according to the definition of $b_i(k)$;

2. At each time $k$, each agent assigns a proper weight to each nonzero $a_{ij}(k)$ such that $\lfloor {L}({k})\rfloor_{ii}$ is no larger than $d_i$ for some constant $d_i<\frac{p^2_i(0)}{4}$;

3. At each time $k$, based on $p_i(k)$, select $p_i(k+1)$ such that $p_i(k+1)\geq b_i(k)$ and $0<p_i(k+1)T<1$. Then calculate $b_i(k+1)$ according to its definition.

Clearly, by selecting proper $p_i(0)$ and nonzero $a_{ij}(k)$, the first two steps can be easily realized.
From the first step, we have that $0<p_i(0)T<1$. Hence from the definition of $e_i(k)$, we have that $0<b_i(0)T<1$ and $b_i(0)\geq p_i(0)$. Then there exists $p_i(1)$ such that $b_i(0)T\leq p_i(1)T<1$. That is, the third step is realized and we have that $0<p_i(0)T\leq b_i(0)T\leq p_i(1)T<1$. By analogy,  for all $i$ and all $k$, the third step can be realized and there exist $p_i(k)$ such that $0<p_i(k)T\leq b_i(k)T\leq p_i(k+1)T<1$. Moreover, from the second step, we have $d_i<\frac{p^2_i(k)}{4}$ for all $i$ and all $k$. That is, Assumption \ref{ass3} is satisfied.

From the design rules above, it can be seen that the gains of the algorithm for all agents need not to be the same or predesigned, and they
can be adjusted by each agent itself based on its own and  neighbors' information.

 {\begin{assumption}\label{ass13}{\rm  Suppose that there exist an infinite time sequence of $k_0,k_1,k_2,\cdots$ and a positive integer $\eta$ such that $k_0=0$, $0<k_{m+1}-k_m\leq \eta$ for all $m$ and the union of the graphs
$\mathcal{G}(k_m),\mathcal{G}(k_m+1),\cdots,\mathcal{G}(k_{m+1}-1)$
has directed spanning trees.}\end{assumption}}

\begin{theorem}\label{theorem2}{\rm Under Assumptions \ref{ass2}-\ref{ass13}, for the multi-agent system (\ref{eq11}) with (\ref{eq526}),  all agents reach a consensus on their position states exponentially fast while their velocities remain in their corresponding constraint sets. Specifically,
  \begin{itemize}
  \item [(a)] there exist a vector $\bar{x}\in \mathbb{R}^r$ and two constants $C>0$ and {$0<\mu\leq1$} such that
 $\|x_i(k)-\bar{x}\|\leq C(1-\mu)^{k}$ for all $i$ and any $k\geq 0$;\footnote{Note that here $\bar{x}$ is the consensus vector for all agents' position states}
\item [(b)] $\lim_{k\rightarrow+\infty}v_i(k)=0$ and $v_i(k)\in V_i$  for all $i,k$.
 \end{itemize}
}\end{theorem}

\begin{remark}\label{remark2}{\rm In \cite{ren1, yangmeng, Lixiangwei, Mengzhaolin, CaoRen}, it is assumed that the constraint  set of each agent is a hypercube and the communication graph is undirected.
It is hard to extend to consider more general nonconvex constraint sets
and directed communication graphs, especially when communication delays are involved. {In addition, in \cite{angelia,LinRen12-a,srivast}, it is assumed that
the states of the agents are constrained in certain convex sets and the dynamics of the agents is in the form of single integrators. Their results cannot be directly applied here as well. The reasons  mainly lie in three aspects. First, the agents in this paper have two different states, position and velocity states, which are not independent but instead coupled in the form of integral. Unlike \cite{angelia,LinRen12-a,srivast}, the position states cannot be directly controlled. Second, the constraint sets are generally nonconvex and are on the velocities. The nonlinearity caused
by the nonconvex constraints on velocities would further lead to a stronger nonlinearity on the position states. Both nonlinearities are different from and more complicated than that caused by convex constraint sets in  \cite{angelia,LinRen12-a,srivast}. Third, the constraint operator adopted in this paper is different from the projection operator in \cite{angelia,LinRen12-a,srivast}. The nonlinear dynamics of the two operators are different in nature.}}\end{remark}

For simplicity of expression, we only discuss the case of $r=1$ in the proof of Theorem \ref{theorem2} and the case of $r>1$ can be analyzed in the same way.

\subsection{Multiple Model Transformations in the Proof of Theorem \ref{theorem2}}

To perform analysis on the closed-loop system (\ref{eq11}) with (\ref{eq526}), we first make multiple model transformations in three steps
so as to use the property of nonnegative matrices to analyze the system stability for all $k\geq0$.

Step 1): From the definition of the constraint operator $\mathrm{S}_{V_i}(\cdot)$, $x$ and $\mathrm{S}_{V_i}(x)$ have the same direction for any nonzero $x$. From (\ref{ei10}), we have
\begin{eqnarray*}\label{lov}\begin{array}{lll}&\mathrm{S}_{V_i}[v_i(k)-p_i(k)v_i(k)T+\pi_i(k)]\\=&e_i(k)[v_i(k)-p_i(k)v_i(k)T]+e_i(k)\pi_i(k).\end{array}\end{eqnarray*}
Recall that $0<e_i(k)\leq 1$. It follows from the definition of $b_i(k)$ that
 when $0<p_i(k)T<1$, $0<b_i(k)T<1$ and $b_i(k)\geq p_i(k)$. Because $1-b_i(k)T=e_i(k)(1-p_i(k)T)$, it follows that
  \begin{eqnarray*}\label{lov}\begin{array}{lll}&\mathrm{S}_{V_i}[v_i(k)-p_i(k)v_i(k)T+\pi_i(k)]\\=&v_i(k)-b_i(k)v_i(k)T+e_i(k)\pi_i(k).\end{array}\end{eqnarray*}
{Let {$\phi(k)=[x_1^T(k), v_1^T(k), \cdots, x_n^T(k), v_n^T(k)]^T$} and {$E(k)=\mathrm{diag}\{e_1(k),\cdots,e_n(k)\}.$} Then  the system (\ref{eq11}) with (\ref{eq526}) can be written as
\begin{eqnarray}\label{eq23a}\begin{array}{lll}
\phi(k+1)&=&\{\tilde{A}(k)-[E(k)L_0(k)\otimes I_2]
\tilde{B}\}\phi(k)\\
&+&\sum_{m=0}^M[E(k)\Phi_m(k)\otimes I_2]
\tilde{B}\phi(k-m),\end{array}\end{eqnarray}
where $\tilde{A}(k)=\mathrm{diag}\{\tilde{A}_1(k),\cdots, \tilde{A}_n(k)\}$ with $\tilde{A}_i(k)=\begin{bmatrix}1&T\\
0&1-b_i(k)T\end{bmatrix}$, $\tilde{B}=I_n\otimes \begin{bmatrix}0&0\\T&0\end{bmatrix}$, $L_0(k)=\mathrm{\overline{diag}}({L}(k))$, and
the $ij$th entry of $\Phi_m(k)$ $(m=0,1,\cdots,M)$ is
either zero or equal to the weight of the edge $(j,i)$ if
$\tau_{ij}=m$.}


Step 2): To partly decouple the integral relationship of the position and velocity of each agent,
we introduce another variable substitution. Let
{$Q(k)=\mathrm{diag}\{Q_1(k),\cdots,Q_n(k)\}$} with $Q_i(k)=\begin{bmatrix}1&0\\1&\frac{2}{b_i(k)}\end{bmatrix}$, {$\xi(k)=Q(k)\phi(k)$},
$A(k)=\mathrm{diag}\{A_1(k),\cdots, A_n(k)\}$ with $A_i(k)=\begin{bmatrix}1-\frac{b_i(k)T}{2}&\frac{b_i(k)T}{2}\\
\frac{b_i(k)T}{2}&1-\frac{b_i(k)T}{2}\end{bmatrix}$, {$B(k)=\mathrm{diag}\{B_1(k),\cdots, B_n(k)\}$ with $B_i(k)=\begin{bmatrix}0&0\\ 0&\frac{2}{b_i(k)} \end{bmatrix}$ and $F=\begin{bmatrix}0&0\\T&0\end{bmatrix}$}. Clearly, $A(k)=Q(k)\tilde{A}(k)Q^{-1}(k)$.

{It follows that
\begin{eqnarray}\label{eq23a}\hspace{-0.2cm}\begin{array}{lll}
&&\xi(k+1)=Q(k+1)\phi(k+1)\\&=&Q(k+1)Q(k)^{-1}Q(k)\{\tilde{A}(k)\\&-&[E(k)L_0(k)\otimes I_2]
\tilde{B}\}Q(k)^{-1}Q(k)\phi(k)\\
&+&Q(k+1)Q^{-1}(k)Q(k)\\&\times&\sum_{m=0}^M[E(k)\Phi_m(k)\otimes I_2]\tilde{B}Q^{-1}(k)Q(k)\phi(k-m)\\&=&{Q(k+1)Q^{-1}(k)\{A(k)}\\&-&{B(k)[E(k)L_0(k)\otimes F]
\}\xi(k)+Q(k+1)Q^{-1}(k)}\\&\times&\sum{_{m=0}^MB(k)[E(k)\Phi_m(k)\otimes F]
\xi(k-m)}.\end{array}\end{eqnarray}}Here, it should be emphasized that {\begin{eqnarray*}&&Q^{-1}(k)B(k)[E(k)\Phi_m(k)\otimes F]
\xi(k-m)\\&=&Q^{-1}(k)B(k)[E(k)\Phi_m(k)\otimes F]
Q(k-m)\phi(k-m)\\&=&Q^{-1}(k)B(k)[E(k)\Phi_m(k)\otimes F]
Q(k)\phi(k-m)\end{eqnarray*}} for all $k\geq0$ and all $0\leq m\leq M$, because the entries of the even columns of $B(k)[E(k)\Phi_m(k)\otimes F]$ are all 0 and the odd entries of $\xi(k-m)$, $\phi(k-m)$ and $Q(k)\phi(k-m)$ are all equal correspondingly.

Step 3) To use the property of nonnegative matrices for the analysis of the system,
we introduce an augmented system.
Similar to \cite{lin74}, let
$Z(k)=[\xi^T(k),\xi^T(k-1),\cdots,\xi^T(k-M)]^T$ and $\Psi(k)$ be a matrix composed of $(M+1)^2$ $2n\times 2n$ square blocks such that
$[\Psi(k)]_{11}=A(k)-B(k)[E(k)L_0(k)\otimes F]+B(k)[E(k)\Phi_0(k)\otimes F]$, $[\Psi(k)]_{1l}=B(k)[E(k)\Phi_l(k)\otimes F]$, $[\Psi(k)]_{l,l-1}=I_{2n}$ for all $2\leq l\leq M+1$ and all other blocks are zero matrices,
 where $[\Psi(k)]_{ij}$ denotes the $(i,j)$th block. It follows that
\begin{eqnarray}\label{e37}
Z(k+1)=\mathrm{diag}\{ Q(k+1)Q^{-1}(k), I_{2nM}\}\Psi(k) Z(k).
\end{eqnarray}


\begin{remark}{\rm Since all introduced transformation matrices are nonsingular, the system (\ref{eq11}) with (\ref{eq526}) is equivalent to the system (\ref{e37}) without information loss. }\end{remark}

\subsection{Consensus Stability Analysis in the Proof of Theorem \ref{theorem2}}

Let {$\Gamma(k,s)=\prod_{i=s}^k\mathrm{diag}\{ Q(i+1)Q^{-1}(i), I_{2nM}\}\Psi(i)$} be the transition matrix of the system (\ref{e37}). The relation between $Z(k+1)$ and $Z(s)$ for all $k\geq s\geq0$ can be described by
\begin{eqnarray}\label{e371}
Z(k+1)=\Gamma(k,s) Z(s).
\end{eqnarray}
In the following, we will perform consensus analysis on the system (\ref{e371}).
Below are
some lemmas prepared for the main theorem. Specifically, Lemma \ref{lemma1a} studies the stochasticity of the matrices $\mathrm{diag}\{Q(k+1)Q^{-1}(k), I_{2nM}\}\Psi(k)$ and $\Gamma(k,s)$, and the lower boundedness of the scaling factor $e_i(k)$ and an auxiliary matrix that will be used for the proof of Lemma \ref{lemma87}.
Lemma \ref{lemma87} proves that there exists at least one column among the first $2n$ columns of $\Gamma(k,s)$ such that each entry
is larger than some positive constant when $k-s$ is sufficiently large. Lemma \ref{lemma88} proves that
 all columns of the transition matrix $\Gamma(k,s)$ exponentially tend to the same as $k\rightarrow+\infty$.

\begin{lemma}\label{lemma1a}{\rm Under Assumptions \ref{ass2} and \ref{ass3},
\begin{itemize}
\item [(a).] $\mathrm{diag}\{Q(k+1)Q^{-1}(k), I_{2nM}\}\Psi(k)$ and $\Gamma(k,s)$ are stochastic matrices for any $k\geq s\geq0$;
\item [(b).] For all $i$, when $v_i(k)-p_i(k)v_i(k)T+\pi_i(k)=S_{V_i}[v_i(k)-p_i(k)v_i(k)T+\pi_i(k)]$, {$e_i(k)=1,$} and
when $[v_i(k)-p_i(k)v_i(k)T+\pi_i(k)]\neq S_{V_i}[v_i(k)-p_i(k)v_i(k)T+\pi_i(k)]$, $e_i(k)\geq\frac{\underline{\rho}_i}{(\frac{1}{T}+2nTd_{i})\max_j\{\|Z_j(0)\|\}};$
\item [(c).] Let $\Theta(k)=\mathrm{diag}\{ \mathrm{\overline{diag}}\{Q(k+1)Q^{-1}(k)\}, I_{2nM}\}\Psi(k)$. Each nonzero entry of
$\Theta(k)$ is uniformly lower bounded by some positive constant.
\end{itemize}}\end{lemma}
\noindent{\textbf{Proof:}} (a). Let $k\geq s\geq 0$. {Note that 
\begin{eqnarray}\label{eq10013}Q_i(k+1)Q_i^{-1}(k)=\begin{bmatrix}1&0\\1-\frac{b_i(k)}{b_i(k+1)}&\frac{b_i(k)}{b_i(k+1)}\end{bmatrix}.\end{eqnarray}}
Clearly, each row sum of the matrices $Q(k+1)Q^{-1}(k)$ is 1. Calculating $\sum_{l=1}^{M+1}[\Psi(k)]_{1l}\mathbf{1}$, from the property of the graph Laplacian, we have $\sum_{l=1}^{M+1}[\Psi(k)]_{1l}\mathbf{1}=A(k)\mathbf{1}=\mathbf{1}$. That is, each row sum of the first $2n$ row of $\Psi(k)$ is 1. Observing the form of $\Psi(k)$, each of its row sums is 1.

 Recall that when $0<p_i(k)T<1$, $b_i(k)\geq p_i(k)$ and $0<b_i(k)T<1$. Under Assumption \ref{ass3}, we have for all
$i$ and all $k$ \begin{eqnarray}\label{ebb1}\frac{1}{T}>b_i(k+1)\geq p_i(k+1)\geq b_i(k)\geq p_i(k)\end{eqnarray} and \begin{eqnarray}\label{ebb2}p^2_i(k)> 4 d_{i}.\end{eqnarray} Thus, $\frac{b_{i}(k)}{b_{i}(k+1)}\leq 1$, $1-\frac{b_i(k)T}{2}>\frac{1}{2}$ and $b_i^2(k)\geq p^2_i(k)>4 d_{i}$ for all $i$ and all $k$.
It follows that all entries of $Q(k+1)Q^{-1}(k)$ and $\Psi(k)$ are nonnegative. Thus,
the matrices
$Q(k+1)Q^{-1}(k)$, $\Psi(k)$, $Q(k+1)Q^{-1}(k)\Psi(k)$ and $\Gamma(k,s)$ for any $k\geq s\geq 0$ are stochastic matrices.

(b). Since
$Z(k)=\Gamma(k-1,0)Z(0)$, each $Z_i(k)$ is a convex combination of
$Z_j(0)$ for all $i,j$ and thus $\|Z_i(k)\|\leq
\max_j\{\|Z_j(0)\|\}$ for all $i$. It follows from the definitions of $x_i(k)$, $v_i(k)$, $Z(k)$ and $\xi(k)$ that
 {$$\|\pi_i(k)\|\leq 2nT\lfloor {L}(k)\rfloor_{ii}\max_j\{\|Z_j(0)\|\}$$} and
 $\|v_i(k)-p_i(k)v_i(k)T\|\leq \|v_i(k)\|=\|\frac{b_i(k)}{2}(Z_{2i}(k)-Z_{2i-1}(k))\|$. Note that $b_i(k)<\frac{1}{T}$ and $\lfloor {L}(k)\rfloor_{ii}\leq d_{i}$ under Assumption \ref{ass3}. It follows that $\|v_i(k)-p_i(k)v_i(k)T\|\leq \frac{1}{T}\max_j\{\|Z_j(0)\|\}$ and hence
 {$$\|v_i(k)-p_i(k)v_i(k)T+\pi_i(k)\|\leq (\frac{1}{T}+2nTd_{i})\max_j\{\|Z_j(0)\|\}.$$}
When $v_i(k)-p_i(k)v_i(k)T+\pi_i(k)=S_{V_i}[v_i(k)-p_i(k)v_i(k)T+\pi_i(k)]$, {$e_i(k)=1.$}
When $[v_i(k)-p_i(k)v_i(k)T+\pi_i(k)]\neq S_{V_i}[v_i(k)-p_i(k)v_i(k)T+\pi_i(k)]$, under Assumption \ref{ass2}, {$$\|\mathrm{S}_{V_i}[v_i(k)-p_i(k)v_i(k)T+\pi_i(k)]\|\geq \underline{\rho}_i.$$}
Thus,
\begin{eqnarray*}\begin{array}{lll}e_i(k)&=\frac{\displaystyle\|\mathrm{S}_{V_i}[v_i(k)-p_i(k)v_i(k)T+\pi_i(k)]\|}{\displaystyle\|v_i(k)-p_i(k)v_i(k)T+\pi_i(k)\|}\\
&\geq \frac{\underline{\rho}_i}{(\frac{1}{T}+2nTd_{i})\max_j\{\|Z_j(0)\|\}}.\end{array}\end{eqnarray*}

(c). From (\ref{ebb1}) and (\ref{ebb2}), we have $b_i(k)\geq p_i(0)$, $\frac{1}{b_i(k+1)}>T$ and $b_i^2(k)\geq p^2_i(0)>4 d_{i}$. Hence,
$p_i(0)T\leq\frac{b_{i}(k)}{b_{i}(k+1)}$ and
$\frac{b_i(k)T}{2}-\frac{2Td_{i}}{b_i(k)}\geq \frac{p_i(0)T}{2}-\frac{2Td_{i}}{p_i(0)}>0$. Note that $a_{ij}(k)\geq \mu_c$.
It follows that each nonzero entry of $\mathrm{\overline{diag}}\{Q(k+1)Q^{-1}(k)\}$, $\Psi(k)$ and hence $\Theta(k)$ is uniformly lower bounded by some positive constant.
\endproof

\begin{remark}{\rm  Lemma \ref{lemma1a} shows that the matrix $\mathrm{diag}\{Q(k+1)Q^{-1}(k), I_{2nM}\}\Psi(k)$ is stochastic. But due to the existence of communication delays, $\mathrm{diag}\{Q(k+1)Q^{-1}(k), I_{2nM}\}\Psi(k)$ has zero diagonal entries. Also, from (\ref{eq10013}), the entry of the stochastic matrix $Q_i(k+1)Q_i^{-1}(k)$, $1-\frac{b_i(k)}{b_i(k+1)}$, might be  arbitrarily close to zero and hence the nonzero entries of the stochastic matrix $\mathrm{diag}\{Q(k+1)Q^{-1}(k), I_{2nM}\}\Psi(k)$ might not be uniformly lower bounded by a positive constant. The role of such nonzero entries might have no difference with the zero entries in the matrix $\Gamma(k,s)$, which might make all columns of $\Gamma(k,s)$ not converge to a common vector. For the convergence analysis of the stochastic matrices, most of the existing results
on delay-related consensus require the number of possible nonzero entries to be finite (e.g., \cite{xiao,lin74,LinRen12-a}) and existing approaches
based on the results in \cite{wolf} require the stochastic matrices to have positive
diagonal entries and their nonzero entries to be uniformly lower bounded by a positive constant.
 Though the results of \cite{NedicOzdaglar10}
allow for an infinite number of edge weights and zero diagonal entries, the union of the communication graphs among each certain time interval is assumed to be strongly connected and each nonzero entry of the stochastic matrices is assumed to be uniformly lower bounded by a positive constant. As a result, existing results cannot be directly applied to deal with the problem studied in this paper. In addition, in our algorithm, it is not required that
the feedback gains for all agents be the same, which distinguishes our results from the works on double-integrator consensus \cite{lin74,hong,xie}, where
the feedback gains
must be uniform among the agents.
}\end{remark}


\begin{lemma}\label{lemma87}{\upshape Under Assumptions \ref{ass2}-\ref{ass13},
 there exists  a positive integer $h\in \{1,\cdots,2n\}$ and a number  $0<\hat{\mu}\leq 1$ such that
$\lfloor\Gamma(k_{m+\mathrm{\hat{n}}}-1,k_{m})\rfloor_{ih}\geq \hat{\mu}$
  for all $k_m\geq0$ and $i$, where $\mathrm{\hat{n}}\geq4n(M+1)$ is a positive integer.
 }\end{lemma}
 \noindent{\textbf{Proof:}} Define $\bar{\Gamma}(k,s)=\prod_{i=s}^k\Theta(i)$, where $\Theta(i)$ has been defined in Lemma \ref{lemma1a}. Obviously,
 ${\Gamma}(k,s)-\bar{\Gamma}(k,s)$ is a nonnegative matrix. So, to prove this lemma, we only need to prove that
  there exist two positive numbers $h\in \{1,\cdots,2n\}$ and $0<\hat{\mu}\leq 1$ such that
$\lfloor\bar{\Gamma}(k_{m+\mathrm{\hat{n}}}-1,k_{m})\rfloor_{ih}\geq\hat{\mu}$
  for all $k_m\geq0$ and $i$.

   Let $\mathcal{\bar{G}}(k,s)$ be the directed graph
 whose edge weight matrix is $\bar{\Gamma}(k,s)$. By calculations, for all $j\in
\{1,\cdots,n\}$, all $i\in\{1,\cdots,M\}$ and $l\in \{1,2,\cdots,2n\}$, $\lfloor\bar{\Gamma}(k,k)\rfloor_{2j,2j-1}=\lfloor\Theta(k)\rfloor_{2j,2j-1}>0$, $\lfloor\bar{\Gamma}(k,k)\rfloor_{2j-1,2j}=\lfloor\Theta(k)\rfloor_{2j-1,2j}>0$ and
$\lfloor\bar{\Gamma}(k,k)\rfloor_{2in+l,2(i-1)n+l}=\lfloor\Theta(k)\rfloor_{2in+l,2(i-1)n+l}>0$, and hence
 nodes
$2j-1$ and $2j$ are strongly connected and there is an edge
from node $2(i-1)n+j$ to node $2in+j$ in $\mathcal{\bar{G}}(k,k)$.
 Also,  note that $\lfloor\bar{\Gamma}(k,k)\rfloor_{j,j}=\lfloor\Theta(k)\rfloor_{j,j}>0$
for all $j\in \{1,2,\cdots,2n\}$ and all $k$ and hence $\lfloor\bar{\Gamma}(k,s)\rfloor_{j,j}>0$, i.e., each node $j$ has an edge to itself in the graph ${\bar{\mathcal{G}}}(k,s)$,
for all $j\in \{1,2,\cdots,2n\}$ and all $k\geq s\geq 0$. It is clear that there is an edge
from node $j$ to node $2in+j$ in $\mathcal{\bar{G}}(k+s,k)$ for all $s\geq M$, all $j\in \{1,2,\cdots,2n\}$ and all $i\in \{1,2,\cdots,M\}$.
Since the union of
$\mathcal{G}(k_{m}),\cdots,\mathcal{G}(k_{{m}+1}-1)$
has at least a directed spanning tree under Assumption \ref{ass13},
the graphs $\mathcal{\bar{G}}(k_{m+1}-1,k_m)$ and hence $\mathcal{\bar{G}}(k_{m+s}-1,k_m)$ have directed spanning trees for all integers $s\geq 1$.

Let $\emptyset\neq\mathrm{ro}(\mathcal{\bar{G}}(k_{m+c}-1,k_m))\subseteq \{1,\cdots, 2n\}$ be the set of the root nodes of $\mathcal{\bar{G}}(k_{m+c}-1,k_{m})$ for any positive integer $c$, and $\mathrm{len}(\mathcal{\bar{G}}(k_{m+c}-1,k_m))$ be the maximum value of the lengths of directed paths from any node of $\mathrm{ro}(\mathcal{\bar{G}}(k_{m+c}-1,k_{m}))$ to any other node without going through the same node twice in $\mathcal{\bar{G}}(k_{m+c}-1,k_m)$ for any positive integer $c$.
Let $\mathrm{rocen}(\mathcal{\bar{G}}(k_{m+c}-1,k_m),\mathcal{\bar{G}}(k_{s+b}-1,k_s))\subseteq \{1,\cdots, 2n\}$ be the set of the nodes of the graph $\mathcal{\bar{G}}(k_{s+b}-1,k_s)$ for two positive integers $s,b$ such that each node has an edge from the set of $\mathrm{ro}(\mathcal{\bar{G}}(k_{m+c}-1,k_m))$ and $\mathrm{ro}(\mathcal{\bar{G}}(k_{m+c}-1,k_m))\cap\mathrm{rocen}(\mathcal{\bar{G}}(k_{m+c}-1,k_m),\mathcal{\bar{G}}(k_{s+b}-1,k_s))=\emptyset$, and  $\emptyset\neq\mathrm{\overline{ro}}(\mathcal{\bar{G}}(k_{m+c}-1,k_m))\subseteq \{1,\cdots, 2n\}$ be the set of the nodes of $\mathcal{\bar{G}}(k_{m+c}-1,k_{m})$ such that each node has an edge from the set of $\mathrm{ro}(\mathcal{\bar{G}}(k_{m+c}-1,k_m))$. Clearly,
$\mathrm{{ro}}(\mathcal{\bar{G}}(k_{m+c}-1,k_m))\subseteq \mathrm{\overline{ro}}(\mathcal{\bar{G}}(k_{m+c}-1,k_m))$.
It should be noted here that the definitions of these sets, only the first $2n$ nodes are considered in the graphs.

Recall that there is an edge
from node $j$ to node $2in+j$ in $\mathcal{\bar{G}}(k+s,k)$ for all $s\geq M$, all $j\in \{1,2,\cdots,2n\}$ and all $i\in \{1,2,\cdots,M\}$.
It is clear that $\mathrm{num}(\mathrm{\overline{ro}}(\mathcal{\bar{G}}(k_{m+2+M}-1,k_{m+1})))\geq 1$ and $\mathrm{len}(\mathcal{\bar{G}}(k_{m+2+M}-1,k_{m+1}))\leq 2n$, where
 $\mathrm{num}(\cdot)$ be the number of the elements of a given set. Note that $\lfloor\bar{\Gamma}(k,k)\rfloor_{j,j}>0$
for all $j\in \{1,2,\cdots,2n\}$ and note that $\mathcal{\bar{G}}(k_{m+s}-1,k_m)$ have directed spanning trees for all integers $s\geq 1$.
From graph theory, if $\mathrm{num}(\mathrm{rocen}(\mathcal{\bar{G}}(k_{m+2+M}-1,k_{m+1}),\mathcal{\bar{G}}(k_{m+3+M}-1,k_{m+2+M})))=0$,
$\mathrm{num}(\mathrm{\overline{ro}}(\mathcal{\bar{G}}(k_{m+3+M}-1,k_{m+1})))\geq 2$ and
$\mathrm{len}(\mathcal{\bar{G}}(k_{m+3+M}-1,k_{m+1}))\leq 2n+1-\mathrm{num}(\mathrm{\overline{ro}}(\mathcal{\bar{G}}(k_{m+2+M}-1,k_{m+1})))\leq 2n$.
If $\mathrm{num}(\mathrm{rocen}(\mathcal{\bar{G}}(k_{m+2+M}-1,k_{m+1}),\mathcal{\bar{G}}(k_{m+3+M}-1,k_{m+2+M})))\neq0$,
$\mathrm{num}(\mathrm{\overline{ro}}(\mathcal{\bar{G}}(k_{m+3+M}-1,k_{m+1})))\geq 2$ and
$\mathrm{len}(\mathcal{\bar{G}}(k_{m+3+M}-1,k_{m+1}))\leq\mathrm{len}(\mathcal{\bar{G}}(k_{m+2+M}-1,k_{m+1}))-1\leq 2n-1$.
If $\mathrm{num}(\mathrm{rocen}(\mathcal{\bar{G}}(k_{m+3+M}-1,k_{m+1}),\mathcal{\bar{G}}(k_{m+4+M}-1,k_{m+3+M})))=0$,
$\mathrm{num}(\mathrm{\overline{ro}}(\mathcal{\bar{G}}(k_{m+4+M}-1,k_{m+1})))\geq 3$ and
$\mathrm{len}(\mathcal{\bar{G}}(k_{m+4+M}-1,k_{m+1}))\leq 2n+1-\mathrm{num}(\mathrm{\overline{ro}}(\mathcal{\bar{G}}(k_{m+3+M}-1,k_{m+1})))\leq 2n-1$. If $\mathrm{num}(\mathrm{rocen}(\mathcal{\bar{G}}(k_{m+3+M}-1,k_{m+1}),\mathcal{\bar{G}}(k_{m+4+M}-1,k_{m+3+M})))\neq 0$,
$\mathrm{num}(\mathrm{\overline{ro}}(\mathcal{\bar{G}}(k_{m+4+M}-1,k_{m+1})))\geq 3$ and
$\mathrm{len}(\mathcal{\bar{G}}(k_{m+4+M}-1,k_{m+1}))\leq \mathrm{len}(\mathcal{\bar{G}}(k_{m+4+M}-1,k_{m+1}))-1\leq 2n-1$.

 By analogy, there exists a positive integer $\bar{n}> 2n+M$ such that $\mathrm{len}(\mathcal{\bar{G}}(k_{m+1+\bar{n}}-1,k_{m+1}))=1$. Note that $\lfloor\bar{\Gamma}(k,k)\rfloor_{j,j}>0$. From graph theory, there exists a positive integer $\hat{n}>4n(M+1)$ such that $\mathrm{len}(\mathcal{\bar{G}}(k_{m+\hat{n}}-1,k_{m}))=1$. That is, $\lfloor\bar{\Gamma}(k_{m+\mathrm{\hat{n}}}-1,k_{m})\rfloor_{ih}> 0$ for some $h\in \{1,\cdots,2n\}$ and all $i$. Note from Lemma \ref{lemma1a} that ${\Gamma}(k,s)$ is a stochastic matrix and each nonzero entry of $\bar{\Gamma}(k_{m+\mathrm{\hat{n}}}-1,k_{m})$ is lower bounded by a positive constant. Also, recall that ${\Gamma}(k,s)-\bar{\Gamma}(k,s)\geq 0$. It follows that
$\lfloor\bar{\Gamma}(k_{m+\mathrm{\hat{n}}}-1,k_{m})\rfloor_{ih}\geq\hat{\mu}$  for all $i$ and some $0<\hat{\mu}\leq 1$.
\endproof


\begin{lemma}\label{lemma88}{\upshape Under Assumptions \ref{ass2}-\ref{ass13},

\noindent(a) there exists a constant $0\leq\rho_i(s)\leq1$ for all $i\in \{1,\cdots,2n(M+1)\}$ and any given $s\geq0$ such that
$\lim\limits_{k\rightarrow+\infty}\lfloor\Gamma(k,s)\rfloor_{hi}=\rho_i(s)$
for all $h$ where $\sum_{i=1}^{2n(M+1)}\rho_i(s)=1$.

\noindent(b) there is a constant {$0<\hat{\mu}\leq1$} such that
$\max_i|\lfloor
\Gamma(k,s)\rfloor_{ij}-\rho_i(s)|\leq C_0(1-\hat{\mu})^{\frac{k}{\mathrm{\hat{n}}\eta}}$ for all $k\geq s$, where {$C_0=(1-\hat{\mu})^{-2}$} and $\mathrm{\hat{n}}\geq4n(M+1)$.}\end{lemma}

\noindent{\textbf{Proof:}} Since the transition matrix $\Gamma(k,s)$ is a stochastic matrix for all $k\geq s\geq 0$ from Lemma \ref{lemma1a}, each entry of the transition matrix
$\Gamma(k,s)$, $\lfloor\Gamma(k,s)\rfloor_{hi}$, is a convex
combination of the entries of the $i$th column of $\Gamma(k-1,s)$.
It follows that
$\max_{h}\{\lfloor\Gamma(k,s)\rfloor_{hi}\}\leq
\max_{h}\{\lfloor\Gamma(k-1,s)\rfloor_{hi}\}$ and $\min_{h}\{\lfloor\Gamma(k-1,s)\rfloor_{hi}\}\leq
\min_{h}\{\lfloor\Gamma(k,s)\rfloor_{hi}\}$.
 Thus, the limits of
$\max_h\{\lfloor\Gamma(k,s)\rfloor_{hi}\}$ and $\min_h\{\lfloor\Gamma(k,s)\rfloor_{hi}\}$ exist, denoted by
$0\leq\theta_i\leq1$ and $0\leq\sigma_i\leq1$ respectively. That is, there exists a positive
number $T_0>s$ for any $\epsilon>0$ such that for
all $k>T_0$
{{\begin{eqnarray}\label{eee1}|\max_{h}\{\lfloor\Gamma(k,s)\rfloor_{hi}\}-\theta_i|<\epsilon\end{eqnarray}}} and {\begin{eqnarray}\label{eee2}|\min_{h}\{\lfloor\Gamma(k,s)\rfloor_{hi}\}-\sigma_i|<\epsilon.\end{eqnarray}}

(a) To prove statement (a), we only need to prove that $\theta_i=\sigma_i$ for all $i$.
Suppose that $\theta_i\neq\sigma_i$ for some integer $i$. From Lemma \ref{lemma87},  there exists a positive integer
 $q\in \{1,\cdots,2n\}$ and a positive number $0<\hat{\mu}\leq 1$ such that
$\lfloor\Gamma(k_{m+\mathrm{\hat{n}}}-1,k_{m})\rfloor_{iq}\geq\hat{\mu}$ for all $k_m\geq 0$ and $i$, where $\mathrm{\hat{n}}\geq4n(M+1)$ is a positive integer. Note that
\begin{eqnarray*}\label{eq09c1}\begin{array}{lll}&&\lfloor
\Gamma(k_{m+\mathrm{\hat{n}}}-1,s)\rfloor_{ij}\\
&=&\sum\limits_{l=1,l\neq q}^{2n(M+1)}\lfloor\Gamma(k_{m+\mathrm{\hat{n}}}-1,k_m)\rfloor_{il}
\lfloor\Gamma(k_m-1,s)\rfloor_{lj}\\&+&\lfloor\Gamma(k_{m+\mathrm{\hat{n}}}-1,k_m)\rfloor_{iq}
\lfloor\Gamma(k_m-1,s)\rfloor_{qj}.
\end{array}\end{eqnarray*}
Take $\epsilon<\frac{\hat{\mu}(\theta_j-\sigma_j)}{2(2-\hat{\mu})}$ and $k_m-1>T_0$.
If $\lfloor\Gamma(k_m-1,s)\rfloor_{qj}\leq \frac{\theta_j+\sigma_j}{2}$, we have
\begin{eqnarray*}\label{eq09c}\begin{array}{lll}&&\lfloor
\Gamma(k_{m+\mathrm{\hat{n}}}-1,s)\rfloor_{ij}\\&=& \sum\limits_{l=1,l\neq q}^{2n(M+1)}\lfloor\Gamma(k_{m+\mathrm{\hat{n}}}-1,k_m)\rfloor_{il}
\lfloor\Gamma(k_m-1,s)\rfloor_{lj}+(\hat{\mu}\\&+&\lfloor\Gamma(k_{m+\mathrm{\hat{n}}}-1,k_m)\rfloor_{iq}-\hat{\mu})
\lfloor\Gamma(k_m-1,s)\rfloor_{qj}\\
&\leq& \sum\limits_{l=1,l\neq q}^{2n(M+1)}\lfloor\Gamma(k_{m+\mathrm{\hat{n}}}-1,k_m)\rfloor_{il}
(\theta_j+\epsilon)+\frac{\hat{\mu}(\theta_j+\sigma_j)}{2}\\&+&(\lfloor\Gamma(k_{m+\mathrm{\hat{n}}}-1,k_m)\rfloor_{iq}-\hat{\mu})
(\theta_j+\epsilon)\\
&\leq& (1-\hat{\mu})(\theta_j+\epsilon)+\frac{\hat{\mu}(\theta_j+\sigma_j)}{2}\\
&<& \theta_j-\epsilon
\end{array}\end{eqnarray*}
for all $j$, where the first inequality has used (\ref{eee1}) and
the second inequality has used the fact that $\sum\limits_{j=1}^{2n(M+1)}\lfloor\Gamma(k_{m+\mathrm{\hat{n}}}-1,k_m)\rfloor_{ij}=1$. This contradicts with (\ref{eee1}).
Similarly, if $\lfloor\Gamma(k_m-1,s)\rfloor_{qj}> \frac{\theta_j+\sigma_j}{2}$,
\begin{eqnarray*}\label{eq09c}\begin{array}{lll}&&\lfloor
\Gamma(k_{m+\mathrm{\hat{n}}}-1,s)\rfloor_{ij}\\&=& \sum\limits_{l=1,l\neq q}^{2n(M+1)}\lfloor\Gamma(k_{m+\mathrm{\hat{n}}}-1,k_m)\rfloor_{il}
\lfloor\Gamma(k_m-1,s)\rfloor_{lj}+(\hat{\mu}\\&+&\lfloor\Gamma(k_{m+\mathrm{\hat{n}}}-1,k_m)\rfloor_{iq}-\hat{\mu})
\lfloor\Gamma(k_m-1,s)\rfloor_{qj}\\
&\geq& (1-\hat{\mu})(\sigma_j-\epsilon)+\frac{\hat{\mu}(\theta_j+\sigma_j)}{2}\\
&>& \sigma_j+\epsilon,
\end{array}\end{eqnarray*}
for all $j$, where the first inequality has used (\ref{eee2}).
 This also yields a contradiction.
Therefore, there exists a constant $0\leq\rho_i(s)\leq1$ for each $i$ such that
$\lim\limits_{k\rightarrow+\infty}\lfloor\Gamma(k,s)\rfloor_{hi}=\rho_i(s)$
for all $h$. Moreover, from the stochasticity of  $\Gamma(k,s)$, we have $\sum_{i=1}^{2n(M+1)}\rho_i(s)=1$.

(b)  From Lemma \ref{lemma87} again,  there exists a positive integer
 $q\in \{1,\cdots,2n\}$ and a positive number $0<\hat{\mu}\leq1$ such that
$\lfloor\Gamma(k_{h\mathrm{\hat{n}}}-1,k_{(h-1)\mathrm{\hat{n}}})\rfloor_{iq}\geq\hat{\mu}$ for all $i$ and  all positive integers $h>0$. Note that
\begin{eqnarray}\label{eq09c1}\hspace{-0.3cm}\begin{array}{lll}&\lfloor
\Gamma(k_{h\mathrm{\hat{n}}}-1,s)\rfloor_{ij}\\
=&\sum\limits_{l=1,l\neq q}^{2n(M+1)}\lfloor\Gamma(k_{h\mathrm{\hat{n}}}-1,k_{(h-1)\mathrm{\hat{n}}})\rfloor_{il}
\lfloor\Gamma(k_{(h-1)\mathrm{\hat{n}}}-1,s)\rfloor_{lj}\\+&\lfloor\Gamma(k_{h\mathrm{\hat{n}}}-1,k_{(h-1)\mathrm{\hat{n}}})\rfloor_{iq}
\lfloor\Gamma(k_{(h-1)\mathrm{\hat{n}}}-1,s)\rfloor_{qj}.
\end{array}\end{eqnarray}
{From (\ref{eq09c1}), we have \begin{eqnarray*}\begin{array}{lll}\max_i\lfloor
\Gamma(k_{h\mathrm{\hat{n}}}-1,s)\rfloor_{ij}\leq (1-\hat{\mu})\max_l \lfloor\Gamma(k_{(h-1)\mathrm{\hat{n}}}-1,s)\rfloor_{lj}\\+\hat{\mu}
\lfloor\Gamma(k_{(h-1)\mathrm{\hat{n}}}-1,s)\rfloor_{qj}\end{array}\end{eqnarray*} and \begin{eqnarray*}\begin{array}{lll}\min_i\lfloor
\Gamma(k_{h\mathrm{\hat{n}}}-1,s)\rfloor_{ij}\geq (1-\hat{\mu})\min_l \lfloor\Gamma(k_{(h-1)\mathrm{\hat{n}}}-1,s)\rfloor_{lj}\\+\hat{\mu}
\lfloor\Gamma(k_{(h-1)\mathrm{\hat{n}}}-1,s)\rfloor_{qj}.\end{array}\end{eqnarray*} Hence, \begin{eqnarray*}\begin{array}{lll}\max_i\lfloor
\Gamma(k_{h\mathrm{\hat{n}}}-1,s)\rfloor_{ij}-\min_i\lfloor
\Gamma(k_{h\mathrm{\hat{n}}}-1,s)\rfloor_{ij}\leq (1-\hat{\mu})\\\times(\max_l \lfloor\Gamma(k_{(h-1)\mathrm{\hat{n}}}-1,s)\rfloor_{lj}-\min_l \lfloor\Gamma(k_{(h-1)\mathrm{\hat{n}}}-1,s)\rfloor_{lj}).\end{array}\end{eqnarray*} Thus, \begin{eqnarray*}\begin{array}{lll}\max_i\lfloor
\Gamma(k_{h\mathrm{\hat{n}}}-1,s)\rfloor_{ij}-\min_i\lfloor
\Gamma(k_{h\mathrm{\hat{n}}}-1,s)\rfloor_{ij}\leq (1-\hat{\mu})^{h-1}\\\times(\max_l \lfloor\Gamma(k_{\mathrm{\hat{n}}}-1,s)\rfloor_{lj}-\min_l \lfloor\Gamma(k_{\mathrm{\hat{n}}}-1,s)\rfloor_{lj}).\end{array}\end{eqnarray*}
Since $\max_l \lfloor\Gamma(k_{\mathrm{\hat{n}}}-1,s)\rfloor_{lj}-\min_l \lfloor\Gamma(k_{\mathrm{\hat{n}}}-1,s)\rfloor_{lj}\leq1$ from the stochasticity and nonnegativity of the matrix $\Gamma(k_{\mathrm{\hat{n}}}-1,s)$, \begin{eqnarray*}\begin{array}{lll}\max_{i}\{\lfloor\Gamma(k,s)\rfloor_{ij}\}\leq
\max_{i}\{\lfloor\Gamma(k-1,s)\rfloor_{ij}\}\end{array}\end{eqnarray*} and \begin{eqnarray*}\begin{array}{lll}\min_{i}\{\lfloor\Gamma(k-1,s)\rfloor_{ij}\}\leq
\min_{i}\{\lfloor\Gamma(k,s)\rfloor_{ij}\}\end{array}\end{eqnarray*} for all $k>s$, we have $\max_i\lfloor
\Gamma(k,s)\rfloor_{ij}-\min_i\lfloor
\Gamma(k,s)\rfloor_{ij}\leq\max_i\lfloor
\Gamma(k_{h\mathrm{\hat{n}}}-1,s)\rfloor_{ij}-\min_i\lfloor
\Gamma(k_{h\mathrm{\hat{n}}}-1,s)\rfloor_{ij}\leq (1-\hat{\mu})^{h-1}$ for all $k_{h\mathrm{\hat{n}}}\leq k\leq k_{(h+1)\mathrm{\hat{n}}}-1$.
 Since $h>\frac{k}{\mathrm{\hat{n}}\eta}-1$ for all $k_{h\mathrm{\hat{n}}}\leq k\leq k_{(h+1)\mathrm{\hat{n}}}-1$,
we have $\max_i|\lfloor
\Gamma(k,s)\rfloor_{ij}-\rho_i(s)|\leq (1-\hat{\mu})^{h-1}\leq C_0(1-\hat{\mu})^{\frac{k}{\mathrm{\hat{n}}\eta}}$ where $C_0=(1-\hat{\mu})^{-2}$.} \endproof

\begin{remark}{\rm From Lemma \ref{lemma88}, it can be seen that each row of the product of stochastic matrices exponentially converge to a certain vector when the union of the edges whose weights are lower bounded by a certain positive constant among each certain time interval has a spanning tree, even when the stochastic matrices have zero
diagonal entries and some of their nonzero entries are arbitrarily close to zero.
}\end{remark}

\noindent\textbf{Proof of Theorem \ref{theorem2}}: Under Assumptions \ref{ass2}-\ref{ass13}, from Lemmas \ref{lemma1a} and \ref{lemma88}, $\Gamma(k,s)$ is a stochastic matrix for any $k\geq s$, and there are constants $0\leq\rho_i(s)\leq1$ and {$0<\hat{\mu}\leq1$} for all $i\in \{1,\cdots,2n(M+1)\}$ such that $\sum_{i=1}^{2n(M+1)}\rho_i(s)=1$ and $\max_i|\lfloor
\Gamma(k,s)\rfloor_{ij}-\rho_i(s)|\leq C_0(1-\hat{\mu})^{\frac{k}{\mathrm{\hat{n}}\eta}}$ for all $k\geq s\geq 0$, where  {$C_0=(1-\hat{\mu})^{-2}$} and $\mathrm{\hat{n}}\geq 4n(M+1)$. Since the initial conditions of $x_i(k)$ and $v_i(k)$ for all $k\leq 0$ satisfy the dynamics of (\ref{eq11}), from (\ref{e371}) and the definitions of $Z(k)$ and $\xi(k)$, the solution of $Z(k)$ exists for any given $k\geq0$. Let $s=M$ and $\bar{x}=\sum_{j=1}^{2(M+1)n}\rho_j(s)Z_j(s)$. Thus, $\|Z_i(k)-\bar{x}\|=\|\sum_{j=1}^{2(M+1)n}\lfloor\Gamma(k-1,s)\rfloor_{ij}Z_j(s)-\sum_{j=1}^{2(M+1)n}\rho_j(s)Z_j(s)\|\leq\sum_{j=1}^{2n(M+1)}\|\lfloor
\Gamma(k-1,s)\rfloor_{ij}-\rho_j(s)\|\|Z_j(s)\|\leq C_0(1-\hat{\mu})^{\frac{k-1}{\mathrm{\hat{n}}\eta}}\sum_{j=1}^{2n(M+1)}\|Z_j(s)\|$. As a result,
$\lim_{k\rightarrow+\infty}\|Z_i(k)-\bar{x}\|=0$ and hence $\lim_{k\rightarrow+\infty}[x_i(k)-\bar{x}]=\lim_{k\rightarrow+\infty}v_i(k)=0$. Also, note that
all $\|Z_i(k)-\bar{x}\|$ are bounded for all $0\leq k\leq s$ and all $i$. Since $0<\hat{\mu}\leq1$, it follows that there exist two constants $C>0$ and $0<\mu\leq 1$ such that
$\|x_i(k)-\bar{x}\|\leq C(1-{\mu})^k$ for any $k\geq 0$ and all $i$.
  \endproof

 \section{A Numerical Example}
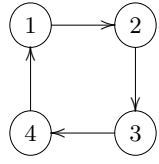
\begin{figure}
\begin{center}
\xymatrix{
    *++[o][F-]{1}\ar[r]
    & *++[o][F-]{2}\ar[d]\\
    *++[o][F-]{4}\ar[u]
    & *++[o][F-]{3}\ar[l]
  }
\end{center}
\vspace{0.2cm}
\caption{One directed graph.} \label{fig:1}
\end{figure}

\begin{figure}
\centering
\includegraphics[width=3.5in]{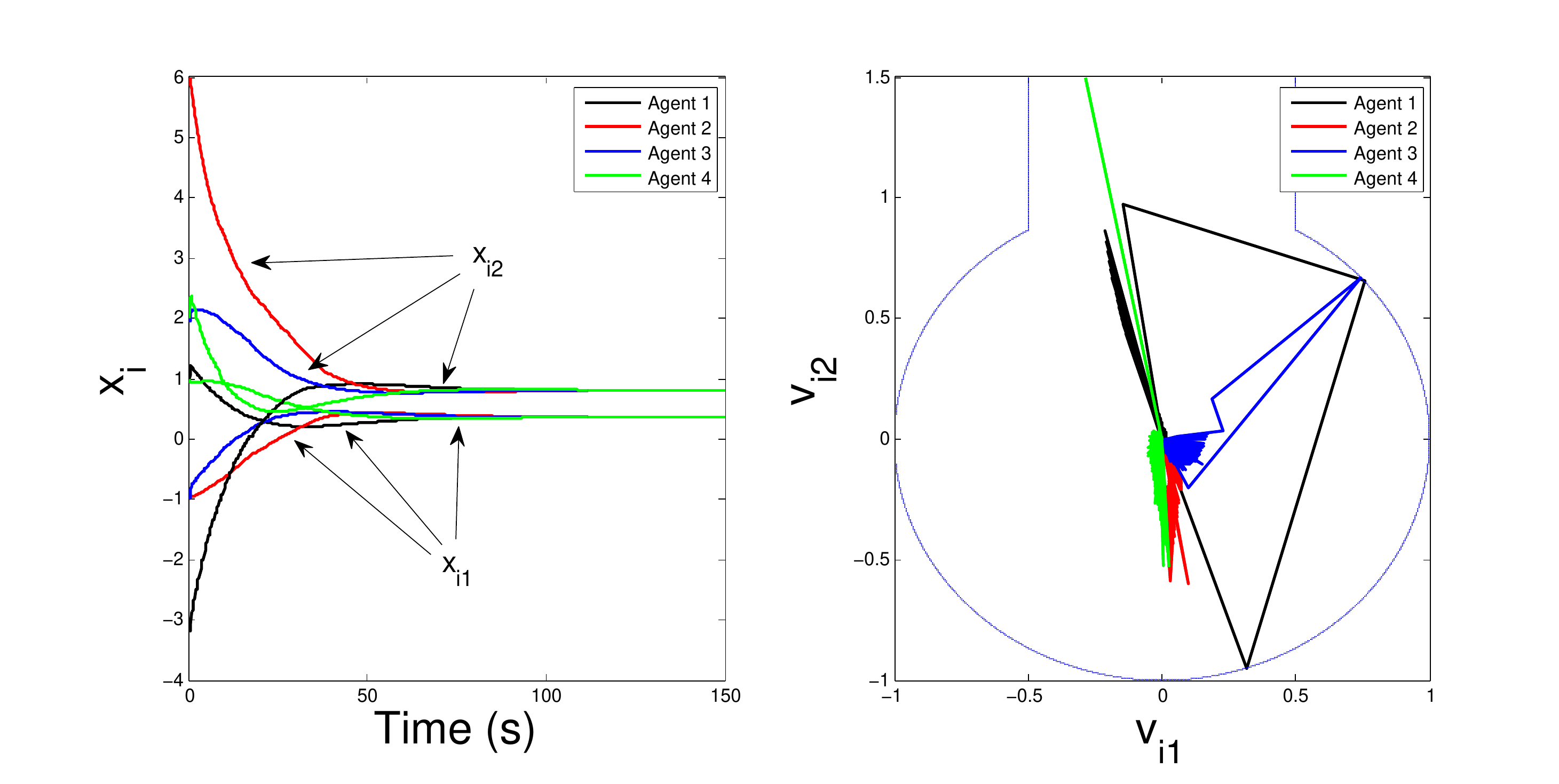} \\
\vspace{-0.4cm}
\caption{Trajectories of all agents}
\label{fig:4}
\end{figure}

Consider a multi-agent system consisting of 4 agents in a plane. The velocity of each agent $v_i$ is constrained to lie in a nonempty nonconvex set $V_i=\{x\mid \|x\|\leq 1\}\cup \{x\mid -0.5\leq [1,0]^Tx\leq 0.5, 0\leq [0,1]^Tx\leq 1.5\}$ for all $i$. At each time, only one edge of the graph shown in Fig. \ref{fig:1} is available to transmit the information and the switching sequence of the edges is (1,2), (2,3), (3,4), (4,1), (1,2), $\cdots$.
 The weight of each edge is 0.5. The
sample time of the system is $T=0.2$ $s$. The initial conditions of all agents are taken as $x_i(k)=x_i(0)$ and $v_i(k)=0$ for all $k<0$. The delay for
the edge (1,2) is $T$ $s$, for the edge (2,3) is $2T$ $s$ and for the edges (3,4) and (4,1) is $3T$ $s$. {According to the design rule of $p_i(k)$ that satisfies Assumption \ref{ass3}, the parameters of the control algorithms
(\ref{eq526}) are taken as $p_i(0)=1.5$ and $p_i(k)=b_i(k-1)$ for all $i$ and all $k\geq 1$.} {Clearly, Assumptions \ref{ass2}-\ref{ass13} are all satisfied.} Fig. \ref{fig:4} shows
the simulation results of the multi-agent system (\ref{eq11}) with (\ref{eq526}). It is clear that all agents eventually reach a consensus while their velocities remain their constraint
sets, which is consistent with Theorem \ref{theorem2}. 

\section{Conclusions}
In this paper, a distributed velocity-constrained consensus problem has been studied
for discrete-time multi-agent systems. Each agent's velocity is constrained to lie in a nonconvex set.
The communication graph considered is directed coupled with arbitrarily bounded communication delays and
can be arbitrarily switching under the condition that the union of the graphs has
directed spanning trees among each certain time interval. A distributed control algorithm
has been proposed by applying a constrained control scheme using only local information.
  {To analyze the velocity-constrained consensus problem, we have first introduced multiple novel model transformations and selected proper
 control parameters to transform the original system into an equivalent system with stochastic matrices. Due to the existence of communication delays and constraints,
the stochastic matrices are state-dependent and have zero diagonal entries, and their nonzero entries might not be uniformly lower bounded by a positive constant. To overcome these coexisting challenges,
with the help of an auxiliary matrix, we have proved that the transition matrix of the equivalent system has at least one column with all positive entries over a certain time interval and used the convexity of a stochastic matrix to show that all rows of the transition matrix tend to the same exponentially.}


\begin{thebibliography}{xx}


\bibitem{saber} R. Olfati-Saber, J. A. Fax and R. M. Murray,``Consensus and cooperation in
networked multi-agent systems," in Proceedings of the IEEE, vol. 95,
no. 1, pp. 215-233, 2007.




\bibitem{xiao} F. Xiao and L. Wang, ``State consensus for multi-agent systems with switching
topologies and time-varying delays", International Journal of
Control, vol. 79, no. 10, pp. 1277-1284, 2006.





\bibitem{xie} G. Xie and L. Wang. ``Consensus control for a class of networks of dynamic agents", International Journal of Nonlinear and Robust Control,
vol. 17, no. 10-11, pp. 941-959, 2007.

\bibitem{hong}Y. Hong, L. Gao, D. Cheng and J. Hu, ``Lyapunov-Based Approach to Multiagent Systems With Switching Jointly Connected Interconnection". IEEE Transactions on Automatic Control, vol. 52, no. 5, pp. 943-948, 2007.

\bibitem{lin74} P. Lin and Y. Jia.
``Consensus of second-order discrete-time multi-agent systems with
nonuniform time-delays and dynamically changing topologies".
Automatica, vol. 45, no. 9, pp. 2154-2158, 2009.

\bibitem{ren1}W. Ren, ``On Consensus Algorithms for Double-integrator
Dynamics," IEEE Transactions on Automatic Control, vol. 53,
no. 6, pp. 1503-1509, 2008.

\bibitem{yangmeng} T. Yang, Z. Meng, D.V. Dimarogonas, and K.H. Johansson,
``Global consensus for discrete-time multi-agent systems with input saturation constraints", Automatica, vol. 50, no. 2, pp. 499-506, 2014.

\bibitem{Lixiangwei}Y. Li,  J. Xiang, W. Wei, ``Consensus problems for linear time-invariant multi-agent systems with saturation constraints",
IET Control Theory $\&$ Applications, vol. 5, no. 6, pp. 823-829, 2011.

\bibitem{Mengzhaolin} Z. Meng, Z. Zhao, Z. Lin, ``On global leader-following consensus of identical
linear dynamic systems subject to actuator saturation", Systems $\&$ Control Letters, vol. 62, no. 2, pp. 132-142, 2013.

\bibitem{CaoRen} Y. Cao and W. Ren, ``Distributed Coordinated Tracking With Reduced Interaction via a Variable Structure Approach", IEEE Transactions on
Automatic Control, vol. 57, no. 1, pp. 33-48, 2012.

\bibitem{angelia} A. Nedi$\acute{\mathrm{c}}$, A. Ozdaglar, P. A. Parrilo, ``Constrained
consensus and optimization in multi-agent networks", IEEE
Transactions on Automatic Control, vol. 55, no. 4, pp.922-938, 2010.

\bibitem{LinRen12-a} P. Lin and W. Ren, ``Constrained Consensus in Unbalanced Networks With Communication Delays", IEEE Transactions on Automatic Control, vol. 59, no. 3, pp. 775-781, 2014.

\bibitem{srivast} K. Srivastava, A. Nedi$\acute{\mathrm{c}}$, ``Distributed asynchronous constrained
stochastic optimization", IEEE Journal of Selected Topics in Signal
Processing, vol.5, no.4, pp.772-790,2011.
\bibitem{NedicOzdaglar10}A. Nedi$\acute{\mathrm{c}}$ and A. Ozdaglar,``Convergence Rate for Consensus with Delays",
Journal of Global Optimization, vol. 47, no. 3, pp. 437-456, 2010.

\bibitem{plinwrenysong}P. Lin, W. Ren and Y. Song, ``Distributed Multi-agent Optimization Subject to Nonidentical Constraints and Communication Delays", Automatica, vol.65, no.3, pp. 120-131, 2016.

\bibitem{wolf} J. Wolfowitz, ``Products of indecomposable,
aperiodic, stochastic matrices," in Proceedings of American
Mathematical Society, vol. 15, pp. 733-736, 1963.
\bibitem{s10} C. Godsil and G. Royle,
Algebraic Graph Theory. New York: Springer-Verlag, 2001.
\bibitem{s11} R. A. Horn and C. R. Johnson,
Matrix Analysis. Cambridge, U.K.: Cambridge Univ. Press, 1987.


\balance





















\end{thebibliography}
\end{document}